\documentclass[12pt,reqno]{amsart}

 \usepackage{mathrsfs}

\usepackage{amsmath, amsthm, amssymb}
\usepackage{amsfonts}
\usepackage[dvips]{epsfig}
\usepackage{graphicx}
\usepackage{caption}
\usepackage{subcaption}
\usepackage[english]{babel}
\usepackage{hyperref}
\usepackage{url}
\usepackage{xcolor}
\usepackage[none]{hyphenat}
\usepackage{endnotes}
\renewcommand{\footnote}[1]{\endnote{#1}}

\usepackage{tikz}
\usepackage{rotating}
\usepackage[utf8]{inputenc}
\usepackage{cite}
\usepackage{amscd}
\usepackage{bm}
\usepackage{enumerate}

\usepackage{verbatim}
\usepackage{hyperref}
\usepackage{amstext}
\usepackage{latexsym}
\theoremstyle{plain}
\newtheorem{theorem}{Theorem}[section]

\newtheorem{corollary}[theorem]{Corollary}

\newtheorem{proposition}[theorem]{Proposition}
\newtheorem{lemma}[theorem]{Lemma}
\newtheorem{remark}[theorem]{Remark}

\numberwithin{theorem}{section} \numberwithin{equation}{section}

\newcommand{\average}{{\mathchoice {\kern1ex\vcenter{\hrule height.4pt
width 6pt depth0pt} \kern-9.7pt} {\kern1ex\vcenter{\hrule height.4pt
width 4.3pt depth0pt} \kern-7pt} {} {} }}

\def\R{\mathbb{R}}

\newcommand{\e }{\varepsilon }

\newcommand{\n }{\nabla }

\renewcommand{\phi}{\varphi}

\renewcommand{\O }{\Omega }

\newcommand{\be}{\begin{equation}}
\newcommand{\ee}{\end{equation}}

\def\proof{\noindent{{\bf Proof. }}}
\def\square{\vbox{
\hrule height .4pt \hbox{\vrule width .4pt height 7pt \kern 7pt
\vrule width .4pt} \hrule height .4pt }}
\def\QED{\hfill {$\square$}\goodbreak \medskip}

\renewcommand{\epsilon}{\varepsilon}

\renewcommand{\u}{{u}}    

\begin{document}

\title[Overdetermined problems  for   fully nonlinear  equations in space forms]
{Overdetermined problems  for   fully nonlinear  equations in space
forms}


\author{Ignace Aristide Minlend}
\address{Department of Quantitative Techniques, Faculty of Economics and Applied Management, University of Douala, BP. 2701 - Douala.}
\email{ignace.a.minlend@aims-senegal.org}




\begin{abstract}
We study  overdetermined problems  for   fully nonlinear   elliptic
equations in subdomains $\O$  of the  Euclidean sphere
$\mathbb{S}^{N}$  and  the hyperbolic space  $\mathbb{H}^{N}$.  We
prove, the existence  of a classical   solution to the underlined
equation forces $\O$ to be a geodesic ball in the ambient space. Our
result extends to fully nonlinear equations, a similar result in the
case of semilinear equations with the Laplace operator due to
Kumaresan and Prajapat.
\end{abstract}
\maketitle


\section{Introduction and main result}
In this  paper,  we are concerned with  the classification  results
for   bounded  and  $C^{2,\alpha}$  subdomains of  the round unit
sphere  $\mathbb{S}^{N}\subset \R^N$ and  the hyperbolic space
$\mathbb{H}^{N}$, admitting solutions for fully nonlinear equations
as
\begin{align}\label{eq:proble1}
  \begin{cases}
 F( \n^2_gu, \nabla_g u) +f(u)=0 & \quad \textrm{ in  }\quad  \Omega\\\vspace{1mm}
u=0&  \quad\textrm{ on  }\quad\partial  \O\\\vspace{1mm} u> 0&
\quad\textrm{ in }\quad \O,\\\vspace{1mm} |\nabla_g u|_g=c_0&
\quad\textrm{ on}\quad  \partial \O.
  \end{cases}
  \end{align}
Here  $g$ denotes the Riemannian metric in $\mathbb{S}^{N}$ or
$\mathbb{H}^{N}$. The function $f$ is locally Lipschitz continuous
on $\R_+$, $c_0\in \R_{+}$ and  $F$ is a function defined on
$S_N\times \R^N$, where  $S_N$ is the space of symmetric $N\times N$
matrices and $F(0,0)=0$.

When the  operator $F$ is replaced  by the Laplace Beltrami
operator,  problem \eqref{eq:proble1}  reads
\begin{align}\label{eq:proble2"}
\begin{cases}
 \Delta_{g} u +f(u)=0 & \quad \textrm{ in  }\quad  \Omega\\\vspace{1mm}
u=0&  \quad\textrm{ on  }\quad\partial  \O\\\vspace{1mm} u> 0&
\quad\textrm{ in }\quad \O,\\\vspace{1mm} |\nabla u|_{g}=c_0&
\quad\textrm{ on}\quad  \partial \O.
  \end{cases}
  \end{align}
and is known as the  Serrin's overdetermined problem in the
classical literature.

In the Euclidean space $\R^N$, Serrin \cite{Serrin}  proved in 1971
that  the existence of a classical solution to  \eqref{eq:proble2"}
implies that  $\O$ is a round ball.   To prove his result Serrin
used the moving planes method  originally  developed  by Alexandrov
\cite{Alexandrov}. In fact Alexandrov introduced the moving plane
method and classified the sphere is the only embedded hypersurface
of $\R^N$ with constant mean curvature. We should stress that soon
after Serrin published his work, Weinberger \cite{Wein} introduced
the so called \textit{$P$-function method} and derived a short proof
of Serrin's result by the use of the maximum principle and
Pohozaev's identity \cite{Poho}. The literature on the moving planes
method is vast. We quote the representative references where the
method was extended in the context of PDE to prove monotonicity and
symmetry properties of positive solutions to some nonlinear second
order
PDEs,\cite{BerestC,DamascelliA,DamascelliB,BGidasW,Gidas,Berestycki}.

For some years  it had been a great  challenge to know whether
Serrin's   result is still valid in many others Riemannian manifolds
instead of  the Euclidean space $\R^N$. This question has been
intensively  studied by several authors in compact manifolds
\cite{FPP, FallIgnace, Sic, VazquezEncicoSalas, DS, Sic,
FallMinlendIWeth2}  and in space  forms in particular \cite{Molzon,
CiraoloLuigi, Ku-Pra}.  In \cite{Ku-Pra},  Kumaresan and Prajapat
proved by the  use of the moving planes method  that if problem
\eqref{eq:proble2"} is solvable in a bounded and regular domain $\O$
of $\mathbb{S}^{N}$ or $\mathbb{H}^{N}$, then $\O$ is a geodesic
ball (with the restriction for the case of the unit sphere  that the
closure  $ \overline{\O}$ is contained in  a  hemisphere). We
emphasize that geodesic balls are not the only Serrin domains in
$\mathbb{S}^{N}$. This was shown in \cite{FallMinlendIWeth2},where
the  we constructed Serrin domains in $\mathbb{S}^{N}$, $N\geq 2$
which bifurcate from symmetric straight tubular neighborhoods of the
equator.

Our goal is to extend  the  result by Kumaresan and Prajapat
\cite{Ku-Pra}  to  the  class of  second order fully nonlinear
operators  $F$ satisfying the following assumptions:

The function $F$ is rotationally invariant, that is

$(\textrm{H}_1)$ $F(\n^2_g(u\circ \phi), \nabla_g (u\circ
\phi))=F(\n^2_gu, \n_g u)\circ \phi$, for  all isometry $\phi$ of
the ambient manifold.

The second assumption  involves the Riemannian Pucci operators and
states that the operator $F$ is uniformly elliptic and Lipschitz
continuous on  $S_N \times \R^N$ in the following sense.

$(\textrm{H}_2)$  there exist numbers $\Lambda \geq \lambda>0$,
$k\geq0$ such that
\begin{equation}\label{eqH2}
\mathcal{P}^{+}_{g, \lambda, \Lambda}(A-B)+k|p-q|_g\geq
F(A,p)-F(B,q)\geq \mathcal{P}^{-}_{g, \lambda,
\Lambda}(A-B)-k|p-q|_g,
\end{equation}
for any $A, B \in S_N$ and any tangent vectors $p$ and $q$.

In \eqref{eqH2},  $|\cdot|_g$ stands for the Riemannian norm,
$\mathcal{P}^{-} _{\lambda, \Lambda}$  and $\mathcal{P}^{+ }
_{\lambda, \Lambda}$  are the   Riemannian Pucci operators   defined
by
\begin{equation}\label{eqpucci}
\mathcal{P}^{+ } _{g, \lambda, \Lambda} (M):= \sup_{A\in
\mathcal{A}^g_{\lambda, \Lambda}} \textrm{tr}(AM) \quad \textrm{and}
\quad \mathcal{P}^{- } _{g, \lambda, \Lambda} (M):= \inf_{A\in
\mathcal{A}^g_{\lambda, \Lambda}} \textrm{tr}(AM),
\end{equation}
where  $\mathcal{A}^g_{\lambda, \Lambda}$ is the set of all
symmetric matrices whose eigenvalues  with respect to $g$ lie to
$[\lambda, \Lambda]$. Moreover, denoting by $\mu_i=\mu^g_i(M)$,
$i=1, \dots, N$  the Riemannian  eigenvalues of $M$,
\begin{equation}\label{eqpucci2}
\mathcal{P}^{-} _{g,\lambda, \Lambda} (M)= \lambda  \sum^N_{\mu_i>0}
\mu_i +\Lambda  \sum^N_{\mu_i<0} \mu_i \quad \textrm{and} \quad
\mathcal{P}^{+} _{g,\lambda, \Lambda} (M)= \Lambda  \sum^N_{\mu_i>0}
\mu_i+\lambda  \sum^N_{\mu_i<0} \mu_i.
\end{equation}
Examples  of operators  satisfying  $(\textrm{H}_1)$ and
$(\textrm{H}_2)$ are  given by
\begin{equation}\label{examples}
F^{\pm} (\n_g^2 u, \n_g u):=\mathcal{P}^{\pm}_{g,\lambda,
\Lambda}(\n_g^2 u)\pm k|\n _gu|_g.
\end{equation}

Problem \eqref{eq:proble1}  was already  considered in  the
Euclidean space $\R^N$ by Silvestre and Sirakov in
\cite{SilvetreBoyan} under  similar assumptions  than
$(\textrm{H}_1)$ and  $(\textrm{H}_2)$. Indeed, it is  proved under
the assumptions  in  \cite{SilvetreBoyan}  that if  there exists a
viscosity solution to  \eqref{eq:proble1}, then $\O$ is a ball and
the solution is  radial.

Our aim is  to use the moving plane method in order to recover the
analogue of this  result  when the ambient space is the  Euclidean
sphere   $\mathbb{S}^{N}$  or the hyperbolic space $\mathbb{H}^{N}$.

The following are our main results.
\begin{theorem}\label{theo1}
Let $\O$ be a bounded  domain in the hyperbolic space
$\mathbb{H}^{N}$.  Assume $(\textrm{H}_1)$, $(\textrm{H}_2)$ and
$F(M, p)$ is continuously  differentiable in $M$. If problem
\eqref{eq:proble1}   has a solution  $u\in
C^{2,\alpha}(\overline{\O})$, then $\O$ is a  geodesic  ball and $u$
is radially symmetric.
\end{theorem}

\begin{theorem}\label{theo1'}
Let $\O\subset \mathbb{S}^{N}$  be  a bounded  domain such that
$\overline{\O}$ is contained in a hemisphere.  Assume
$(\textrm{H}_1)$,  $(\textrm{H}_2)$ and  $F(M, p)$ is continuously
differentiable in $M$. If problem \eqref{eq:proble1}  has  a
solution   $u\in C^{2,\alpha}(\overline{\O})$. Then $\O$ is a
geodesic  ball and $u$ is radially symmetric.
\end{theorem}

The proof of Theorem  \ref{theo1}  is inspired by the argument  by
Birindeli and Demengel \cite{Birindeli-Demengel}, which was  also
adapted in  \cite{SilvetreBoyan}. In  \cite{Birindeli-Demengel}, the
authors proved a Serrin-type symmetry result for operators of the
form $|\n u|^{\alpha} \mathcal{M}^{}_{\lambda, \Lambda}(\n^2 u),$
provided the Pucci operator is  a  small perturbation of the
Laplacian.  In the next theorem, we  obtain  for operators $F$ not
necessarily $C^1$, a similar and intermediate  result which will be
applied  to prove Theorem \ref{theo1}.

\begin{theorem} \label{theo2}
Assume $(\textrm{H}_1)$ and  $(\textrm{H}_2)$.   Assume also that
$u$ is a $C^{2,\alpha}(\overline{\O})$  solution of
\eqref{eq:proble1}.  There exists a  positive constant
$\e_0=\e_0(\alpha)<1$ only depending on $\alpha$ such that if $\e\in
(0, \e_0)$ and   $\Lambda=\lambda(1 +\e)$, then   $\O$ is a geodesic
ball and $u$ is radial.
\end{theorem}

We emphasize that the conclusion of  Theorem  \ref{theo2} holds
under the assumption  that the ratio  $\Lambda/\lambda$ is
sufficiently close to one. We also  remark from the proof of Lemma
\ref{lemma order homoge} in  \cite[Lemma 4.4]{SilvetreBoyan} that if
we replace $\Lambda=\lambda(1 +\e)$ with  $\Lambda=\lambda +\e$, we
still recover  the conclusion of Theorem  \ref{theo2} provided the
ellipticity constants
$\lambda$ and $\Lambda$ are sufficiently close to each other. \\

Theorem   \ref{theo1}  and  Theorem \ref{theo1'}  are stated under
the crucial   assumption  that the operator $F(M,p)$ is $C^1$ in
$M$, and hence does not  include the Pucci equations. Nevertheless,
we  deduce  from   Theorem  \ref{theo2}, the  following   corollary
which  is valid in both  $\mathbb{S}^{N}$  and  $\mathbb{H}^{N}$,
with the restriction  that the closure  $ \overline{\O}$ of the
subdomain $\O$ in  $\mathbb{S}^{N}$   is contained in a  hemisphere.

\begin{corollary} \label{Coro1}
Assume there exists a nonnegative solution $u\in
C^{2,\alpha}(\overline{\O})$ of
\begin{align}\label{eq:prob}
 \begin{cases}
\mathcal{P}^{\pm} _{\lambda, \lambda(1 +\e)} (\n^2 u)+f(u)=0 & \quad
\textrm{ in  }\quad  \Omega\\\vspace{1mm} u=0&  \quad\textrm{ on
}\quad\partial  \O\\\vspace{1mm} |\nabla u|=c_0&  \quad\textrm{
on}\quad  \partial \O.
  \end{cases}
  \end{align}
There exists a positive constant   $\e_0<1$ only depending on
$\alpha$ such that if  $\e\in(0, \e_0)$,  then  $\O$ is a geodesic
ball  and $u$ is radial.
\end{corollary}

We now explain in details our argument  to prove the above results.
As already mentioned previously, Theorem  \ref{theo1} will be deduced from
Theorem  \ref{theo2}.  For the proof of Theorem \ref{theo1'}, we will  transform
problem \eqref{eq:proble1} to an equivalent problem on $\R^N\setminus\{0\}$ via the
stereographic projection  and apply the same argument in the proof of
 Theorem \ref{theo1}, see Section  \ref{prooftheo1'}. \\

To prove Theorem \ref{theo2}, we use the moving planes method which
we describe in details in Section  \ref{Movingplane}. The idea
consists in comparing  the value of the solution of
\eqref{eq:proble1} between two points, one of these points is the
reflection of the other with respect a  "plane". We move the plane
until it reaches a  critical  position which determines the symmetry
of the domain $\O$.  We will apply this method to an equation
involving the Pucci operator $\mathcal{M}^{- } _{\lambda, \Lambda}$,
and  which  we derive   after exploiting the hypothesis
$(\textrm{H}_2)$.  Compared to the Euclidean case, we  are  going to
"move" instead of planes, but complete and totally geodesic
hypersurfaces since the isometry groups of sphere and hyperbolic
space  are both  generated by reflections with respect  to these
hypersurfaces.

When applying the moving planes  method,  two situations often
occur:  either the reflection (with respect to the plane) of a
suitable portion of the domain $\O$ is internally  tangent to the
boundary  $\partial \O$ at some point $P$, or "the plane" is
orthogonal to $\partial \O$ at some point $Q$ (called the corner
point). We handle the first case by applying  Hopf's lemma.  The
second case is more involved and requires a different analysis.
Indeed a  usual tool to deal with the corner point $Q$ is the so
called   "Serrin corner lemma". Unfortunately, the application of
this  lemma fails  in  general for nonlinear equations  as it can be
seen in  \cite{ArmstrongSirakov}  and   \cite[Sec. 4]{SilvetreBoyan}
for the Pucci operators.  Instead of Serrin corner lemma, we use the
fact that  classical solutions to fully nonlinear uniformly elliptic
equations  are $C^{2,\alpha}$ at the boundary and enjoy a Taylor
expansion of order $2+\alpha$, see    \cite[Proposition
2.2]{SilvetreBoyan} and \cite{SilvetreBoyan2}.   The proof of
Theorem \ref{theo2} is achieved applying Proposition \ref{hopfgene}
below, which provides a  non degeneracy result of order strictly
less than $2+\alpha$ under the hypothesis that the  ratio
$\Lambda/\lambda$ is sufficiently  close  to one.

The paper is organized as follows. In Section  \ref{eqpreli}, we set
the preliminaries and give the expression of  the Hessian in the
hyperbolic space $ \mathbb{H}^{N}$. We also provide in Lemma
\ref{Comparapucci}, an inequality result between   the  Pucci
operators  of the hyperbolic space and the Euclidean space $\R^N$.
This inequality will permit us  to run the moving planes method in
Section \ref{Movingplane}  with an  equation only involving the
Pucci operator of $\R^N$.  In section \ref{prooftheo2}, we rule out
the corner situation  via a contradiction argument and hence proving
Theorem \ref{theo2}. The proof of Theorem \ref{theo1} is given in
Section \ref{prooftheo1}, where we show how the $C^1$ assumption
made on  $F$  allows to  choose the  ratio $\Lambda/\lambda$  close
to one. In Section \ref{prooftheo1'}, we  explain  how  we  deduce
the proof of Theorem  \ref{theo1'} from the argument  on the
hyperbolic space.

\bigskip
\noindent \textbf{Acknowledgements}: This work was initiated when
the author won the inaugural  Abbas Bahri Excellence Fellowship at
Rutgers University of New-Jesey, USA. He is greatful to  the
Department of Mathematics for the hospitality and  wishes to thank
his mentors Pr. Yan Yan Li,  Pr. Zheng-chao Han, and Pr. Denis
Kriventsov  for their  helpful suggestions  and insights  throughout
the writing of this paper.

\section{Prelimanaries and notations}\label{eqpreli}

In the following, we  denote by  $ \mathcal{M}_{\lambda,
\Lambda}^{\pm} $  the  Pucci operators  with respect  to the
Euclidian metric in $\R^N$. Similarly for a  $C^2$ function $u$,
$\nabla^2_g u$ and $\nabla_g u$ will denote the  Riemannian Hessian
and the  gradient,  and we will write  $\nabla^2 u$ and $\nabla u$
when the Riemannian metric is the Euclidean one.  If there is no
ambiguity,  we will omit  the subscript  $g$ and simply  write  $
\mathcal{P}_{\lambda, \Lambda}^{\pm} $  to denote  the Riemannian
Pucci operators.

The Riemannian  Hessian of a smooth function $u$ is a $(0, 2)$
tensor   defined by
\begin{equation}\label{RiemHess}
\nabla_g^2 u (X,Y) =X Y(u)-(\n_X Y)(u)
 \end{equation}
 for any vectors fields $X$ and $Y$.

In   local coordinates  we have
$$\nabla_g^2 u=\sum^N_{i,j=1}\biggl(\frac{\partial^2 u}{\partial x_i\partial x_j}-\sum^N_{k=1}\frac{\partial u}{\partial x_k}\Gamma^{k}_{ij}\biggl)dx_i\otimes dx_j,$$
where  $ \Gamma^{k}_{ij}$ are the Christoffel symbols defined by
\begin{equation}\label{ChristSymhy0}
 \Gamma^{k}_{ij}=\frac{1}{2}\sum^N_{r=1}g^{kr}\biggl(\partial_i g_{rj}+ \partial_j g_{ir}-\partial_r g_{ij} \biggl).
 \end{equation}

An eigenvalue of the Hessian $\nabla_g^2 u $ at a point $x$ of a
Riemannian manifold $(\mathcal{M}, g)$ is a number $\mu$  such that
there exists a tangent vector $X\ne0$ and  for any tangent vector
$Y$ at $x$,
 \begin{equation}\label{eigenHess}
\nabla_g^2 u (X,Y)=\mu \langle X, Y\rangle_g.
 \end{equation}
Observe that if $g$ and $h$ are two Riemannian metrics  such that
$g(x)=m(x)h(x)$ for some positive function $m$, denoting by $\mu_g$
and $ \mu _h$ the eigenvalues of the Hessian with respect to $g$ and
$h$ respectively, we have
\begin{equation}\label{relagen}
\mu _h(x)=m(x)\mu_g(x).
\end{equation}

In  what  follows, we  consider the half-space model
$(\mathbb{H}^{N}, g)$ of the hyperbolic space, where
$$\mathbb{H}^{N}= \{(x_1, \dots, x_N)\in \R^N: x_N> 0\}$$   and  $g$
is the metric defined by
\begin{equation}\label{metrichyypp}
g_{ij}=x_N^{-2}\delta_{ij}.
\end{equation}

The Laplace-Beltrami operator on  $(\mathbb{H}^{N}, g)$  is given by
$$\Delta=x_N^2\biggl(\sum_{i=1}^N \frac{\partial^2}{\partial
x_{i}^2}\biggl) +(2-N)x_N\frac{\partial}{\partial x_N}$$
 and from  \eqref{eqpucci2}, \eqref{relagen}  and \eqref{metrichyypp} we have the relation
\begin{equation}\label{relPu}
\mathcal{P}_{\lambda, \Lambda}^{\pm}(\nabla_g^2 u)= x_N^{2} \mathcal{M}_{\lambda, \Lambda}^{\pm}(\nabla_g^2 u).\\
 \end{equation}

Next we provide the expression of the Hessian.

By a  straightforward computation,
\begin{equation}\label{ChristSymhy1}
\Gamma^{k}_{ij}=-(x_N)^{-1}\{\delta_{iN}\delta_{kj}+\delta_{jN}\delta_{ik}-\delta_{kN}\delta_{ij}\},
\end{equation}
from which we deduce
\begin{equation}\label{ChristSymhy}
(\nabla^2_g u)_{ij}=\frac{\partial^2 u}{\partial x_i\partial x_j}+
(x_N)^{-1}\{\frac{\partial u}{\partial
x_j}\delta_{iN}+\frac{\partial u}{\partial
x_i}\delta_{jN}-\frac{\partial u}{\partial x_N}\delta_{ij}\}.
\end{equation}
Hence  \begin{equation}\label{Hessexpress} \nabla^2_g u=\nabla^2
u+x_N^{-1} \mathcal{K}(u),
\end{equation}
where
\begin{equation}\label{inte}
( \mathcal{K}(u))_{ij}:=\frac{\partial u}{\partial
x_j}\delta_{iN}+\frac{\partial u}{\partial
x_i}\delta_{jN}-\frac{\partial u}{\partial x_N}\delta_{ij}=
\begin{cases}
-\dfrac{\partial u}{\partial x_N}\delta _{ij} \quad \textrm{ if }  i\ne N \quad \textrm{and} \quad   j\ne N\vspace{3mm}\\
\dfrac{\partial u}{\partial x_j} \quad \textrm{ if }  i= N \quad \textrm{and} \quad   j\ne N\vspace{3mm}\\
\dfrac{\partial u}{\partial x_i} \quad \textrm{ if }  i\ne N \quad \textrm{and} \quad   j=N\vspace{3mm}\\
\dfrac{\partial u}{\partial x_N} \quad \textrm{ if }  i=N \quad \textrm{and} \quad   j=N. \\
  \end{cases}
\end{equation}

In the following lemma, we give inequality relating the Pucci
operators $ \mathcal{P}_{\lambda, \Lambda}^{\pm}$ and $
\mathcal{M}_{\lambda, \Lambda}^{\pm}.$
\begin{lemma}\label{Comparapucci}
For any   function $u\in C^{2,\alpha}(\overline{\O})$,
\vspace*{0.1cm}
\begin{equation}\label{PucciInegality1}
 \mathcal{P}_{\lambda, \Lambda}^{-}(\nabla^2_g u)-k|\nabla_g u|_g\geq x_N^{2} \mathcal{M}_{\lambda, \Lambda}^{-}(\nabla^2 u) -\mu x_N|\nabla u|
\end{equation}
and
\begin{equation}\label{PucciInegality2}
\mathcal{P}_{\lambda, \Lambda}^{+}(\nabla^2_g u)+k|\nabla_g u|_g\leq
x_N^{2} \mathcal{M}_{\lambda, \Lambda}^{+}(\nabla^2 u) +\mu
x_N|\nabla u|,
\end{equation}
where $\mu:=\mu( \Lambda, N, k)=\Lambda (N-1)+k.$
\end{lemma}
\vspace*{0.1cm}

Though  we  will not use it in this paper, the  inequality
\eqref{PucciInegality2}  is the one  we should  use, in case instead
of  working with the right  hand side,  we choose to the left side
of $(\textrm{H}_2)$.

\proof We have
 \begin{equation}\label{eqpuccirelateR}
\mathcal{P}_{\lambda, \Lambda}^{-}(\nabla^2_g u) =x_N^{2}
\mathcal{M}_{\lambda, \Lambda}^{-}(\nabla^2_g u),\quad \nabla_g u
=x_N^{2} \nabla u\quad\textrm{and}\quad  |\nabla_g u|_g =x_N |\nabla
u|.
\end{equation}

Next, we study  the eigenvalues of   $\mathcal{K}(u)$ in
\eqref{inte}.

Let us  consider  the $N\times N$  matrix $A$
 defined  by
\begin{equation}\label{matrix1}
A=
 \begin{pmatrix}
\delta &  0  & \ldots & \ldots & 0 & a_1\\
0  &  \ddots & \ddots & \ddots & \vdots& \vdots\\
\vdots&  \ddots & \ddots & \ddots & \vdots& \vdots\\
\vdots & \ddots & \ddots & \ddots & 0& \vdots\\
0  &  \ldots & \ldots & 0& \delta& a_{N-1}\\
a_1  & \ldots & \ldots &  \ldots &a_{N-1} & \beta
\end{pmatrix}.
\end{equation}
For   all $k=0,1,...,N-2$,  we introduce the  $N-k$ dimensional
matrix $M_k$ by
\begin{equation}\label{matrix1}
M_k=
 \begin{pmatrix}
\delta &  0  & \ldots & \ldots & 0 & a_1\\
0  &  \ddots & \ddots & \ddots & \vdots& \vdots\\
\vdots&  \ddots & \ddots & \ddots & \vdots& \vdots\\
\vdots & \ddots & \ddots & \ddots & 0& \vdots\\
0  &  \ldots & \ldots & 0& \delta& a_{N-1-k}\\
a_1  & \ldots & \ldots &  \ldots &a_{N-1-k} & \beta
\end{pmatrix}
\end{equation}
so that $M_0=A$.

We have in particular
\begin{equation}\label{eq:pre}
\det(M_{N-2})=\delta\beta-a^{2}_1\quad\textrm{and}\quad\det(M_{N-3})=\delta^2\beta-\delta(a^{2}_1+a^{2}_2).
\end{equation}
The computation of  $\det(M_{k})$ by developing with respect to
column $N-1-k$ yields
\begin{equation}\label{eq:pre1}
\det(M_{k})=\delta\det(M_{k+1})-\delta^{N-k-2}a^{2}_{N-1-k}.
\end{equation}
\textbf{Claim:} We claim that
\begin{equation}\label{eq:claim}
\det(M_{N-1-\ell})=\delta^{\ell}\beta-\delta^{\ell-1}(a^2_1+...+a^2_\ell),\quad
\ell=1,...,N-1.
\end{equation}
Indeed, for $\ell=1$ and $\ell=2$, \eqref{eq:claim} is true by
\eqref{eq:pre}. It is not hard  to show  by recurrence   using
\eqref{eq:pre1} that  the equation \eqref{eq:claim} holds true for
all $\ell=1,...,N-1$.
Equation \eqref{eq:claim} allows us to deduce that\\
\begin{equation}\label{eq:deterfi}
\det(A)=\det(M_{0})=\delta^{N-1}\beta-\delta^{N-2}(a^2_1+...+a^2_{N-1})=\delta^{N-2}[\delta\beta-(a^2_1+...+a^2_{N-1})]
\end{equation}
Now  computing  \eqref{eq:deterfi} with $\delta=-\partial_N
u-\gamma$, $\beta=\partial_N u-\gamma$ and $a_i=\partial_i u$, $i=1,
\dots, N-1$, we get

\begin{equation}\label{eq:detK}
\det(\mathcal{K}(u)-\gamma I_N)=(-1)^N(\gamma+\partial_N
u)^{N-2}(\gamma^{2} -|\n u|^2)
\end{equation}
and we see that the eigenvalues of $\mathcal{K}(u)$ are

\begin{align}\label{eq:eigentK}
\mu_i=-\partial_N u,  \quad i=1, \dots, N-2, \quad \mu_{N-1}=-|\n
u|, \quad \mu_{N}=|\n u|.
\end{align}

At any point $x \in \O$, we have
\begin{equation}\label{PucciK}
\mathcal{M}_{\lambda,
\Lambda}^{-}(\mathcal{K}(u))=\lambda\mu_{N}+\Lambda\mu_{N-1}+
\begin{cases}
-\lambda  (N-2)\partial_N u \quad \textrm{ if }   \partial_N u (x)<0 \vspace{3mm}\\
 -\Lambda  (N-2)\partial_N u \quad \textrm{ if }   \partial_N u (x)> 0.
  \end{cases}
\end{equation}

Using  the fact that  $\partial_N u\leq |\partial_N u|\leq |\n u|$,
$\mu_{N}\geq0$ and $\lambda\leq  \Lambda$, it follows from
\eqref{PucciK} that
\begin{equation}\label{PucciKboud1}
\mathcal{M}_{\lambda,
\Lambda}^{-}(\mathcal{K}(u))\geq\Lambda[\mu_{N-1}-(N-2)|\n u|]=
-\Lambda(N-1)|\n u|.
\end{equation}
Recalling   $$\mathcal{M}_{\lambda, \Lambda}^{-}(A+B)\geq
\mathcal{M}_{\lambda, \Lambda}^{-}(A)+\mathcal{M}_{\lambda,
\Lambda}^{-}(B)$$ for all  symmetric matrices $A, B$, we get from
\eqref{PucciKboud1},  \eqref{Hessexpress} and \eqref{relPu}
\begin{align}\label{Puccicompar}
 \mathcal{P}_{\lambda, \Lambda}^{-}(\nabla^2_g u)-k|\nabla_g u|_g&=x_N^{2} \mathcal{M}_{\lambda, \Lambda}^{-}(\nabla^2_g u)-kx_N |\nabla u|\nonumber\\
 &\geq x_N^{2} \mathcal{M}_{\lambda, \Lambda}^{-}(\nabla^2 u)+ x_N \mathcal{M}_{\lambda, \Lambda}^{-}(\mathcal{K}(u))-kx_N |\nabla u|\nonumber\\
 &\geq  x_N^{2} \mathcal{M}_{\lambda, \Lambda}^{-}(\nabla^2 u) -x_N[\Lambda (N-1)+k]|\nabla u|\nonumber
\end{align}
and hence
\begin{equation}\label{PucciInegality}
\mathcal{P}_{\lambda, \Lambda}^{-}(\nabla^2_g u)-k|\nabla_g u|_g\geq
x_N^{2} \mathcal{M}_{\lambda, \Lambda}^{-}(\nabla^2 u) -x_N\{\Lambda
(N-1)+k\}|\nabla u|,
\end{equation}
which is exactly the inequality \eqref{PucciInegality1}.

It is also plain that
\begin{equation}\label{PucciKb2}
\mathcal{M}_{\lambda,
\Lambda}^{+}(\mathcal{K}(u))\leq[\Lambda(N-1)-\lambda]|\n u|\leq
\Lambda(N-1)|\n u|.
\end{equation}
Now using  $$\mathcal{M}_{\lambda, \Lambda}^{+}(A+B)
\leq\mathcal{M}_{\lambda, \Lambda}^{+}(A)+\mathcal{M}_{\lambda,
\Lambda}^{+}(B),$$ we   get from  \eqref{PucciKb2},
\eqref{Hessexpress} and \eqref{relPu},
\begin{equation}\label{PucciInegality3}
\mathcal{P}_{\lambda, \Lambda}^{+}(\nabla^2_g u)+k|\nabla_g u|_g\leq
x_N^{2} \mathcal{M}_{\lambda, \Lambda}^{+}(\nabla^2 u) +x_N\{\Lambda
(N-1)+k\}|\nabla u|,
\end{equation}

\QED

\section{The moving plane method}\label{Movingplane}
As explained in the introduction,  along the moving "planes" method,
we are  considering reflections  with respect to complete and
totally geodesic hypersurfaces, in contrast to the Euclidean case,
where the reflections are made with respect to  planes.  We then
recall the definition of a  reflection in the setting of  hyperbolic
space, see also  \cite{Ku-Pra}.

Let  $x\in \mathbb{H}^{N}$  and  $\Gamma$  a  complete and totally
geodesic hypersurface of  $\mathbb{H}^{N}$, there exists a unique
point $p\in \Gamma$  such that  $\textrm{d}(x,  p)=\textrm{d}(x,
\Gamma)$. Let $\gamma$  be a geodesic such that  $\gamma(0)=p$  and
$\gamma(t)=x$. The hyperbolic space   $\mathbb{H}^{N}$  being
complete,  the geodesic $\gamma$  can be extended to the whole  $\R$
and we can define the mapping $R_\Gamma: \mathbb{H}^{N}\rightarrow
\mathbb{H}^{N}$  by  $$R_\Gamma(x)=\gamma(-t).$$ The map  $R_\Gamma$
is by definition the reflection with respect to the complete and
totally
geodesic hypersurface  $\Gamma$.\\

Before we properly start the moving "planes" method, we  also make
the following important remark.

Observe that the hyperplanes orthogonal to  $e_1=(1,0,\cdots,0)$,
that is  $T_s:= \{x\in  \mathbb{H}^{N},\quad  x_1=s\}$   are
complete  and totally geodesic hypersurfaces of  $\mathbb{H}^{N}$.

We prove in the first step  that for the fixed direction $e_1$,
there exists a hyperplane $T_{s_*}$  orthogonal to $e_1$ such that
 \begin{equation}\label{eq: syme1}
 u(x)=u(R_{s_*}(x))\quad  \textrm{ for all}\quad   x\in \O,
 \end{equation}
 where $R_s$ denotes the reflection with respect to the hyperplane  $T_s.$\\

Assume  that  \eqref{eq: syme1} is already proved and   consider  a
direction  $e$ different  from   $e_1$. Without loss of generality,
we can assume the $e$ is a unit vector. From the transitive action
of the isometry group of the hyperbolic space \cite[Proposition
1.2.1]{JPau}, there exits an isometry   $\varphi$ of  $
\mathbb{H}^{N}$ transforming  $e$ to a direction $ \widetilde{e}$
parallel to $e_1$.

We define $$\widetilde{u}:=u \circ \varphi^{-1}\quad
\textrm{on}\quad \widetilde{\O}:=\varphi(\O).$$ Then by
$(\textrm{H}_1)$ the function  $\widetilde{u}$ satisfies
$$F( \n^2_g\tilde{u}, \nabla_g \tilde{u})=F( \n^2_g\u, \nabla_g u)\circ \varphi^{-1}=
-f(u)\circ \varphi^{-1}=-f(u\circ \varphi^{-1})=-f(\widetilde{u}) \quad \textrm{ in  }\quad  \widetilde{\O}$$
and we have $\widetilde{u}=0$ on $\partial \widetilde{\O}$ and $\widetilde{u}>0$ in $\widetilde{\O}$.

Next for a point  $p\in \partial\O$, let  $\nu(p)$  be the  normal
vector  at $p$. Since $\varphi$ is an isometry, then the vector
field $$\tilde{\nu}(p):=D\varphi (p)\cdot \nu(p)$$ is the normal
vector to $\partial\tilde{\O}$ at $\varphi (p)$ and we have
$$\frac{\partial \widetilde{u}}{\partial\widetilde{\nu}}(\varphi (p))=
D\widetilde{u}(\varphi (p))\cdot\widetilde{\nu}(p)=Du (p)\circ D\varphi^{-1}(\varphi (p))[(D\varphi(p)\cdot\nu(p)]=\frac{\partial u}{\partial \nu}(p)$$
and we deduce  that $\widetilde{u}$ solves  problem \eqref{eq:proble1} in $\widetilde{\O}$. \\

Hence from the first step, there exists a hyperplane  $T$
orthogonal to  $ \widetilde{e}$  such that
$$\widetilde{u}(x)=\widetilde{u}(R_T x),\quad \textrm{for all  $x\in
\widetilde{\O}$}.$$ This is equivalent to say that
$$
u (y)=u(\varphi^{-1} R_T \varphi (y))\quad \textrm{for all }\quad
y\in \O.
$$
We put  $\Gamma:=\varphi^{-1}(T)$. Then  $\Gamma$ is a complete  and
totally geodesic hypersurface  orthogonal  to $e$. By \cite[Theorem
3]{Ku-Pra}, the reflection with respect to $\Gamma$ is given by
$$R_\Gamma=\varphi^{-1} \circ R_T \circ \varphi$$  and hence $$u
(y)=u( R_\Gamma(y))\quad \textrm{for all }\quad  y\in \O.$$ To
summarize,  we have shown  that for all direction $\gamma$ in
$\R^N$, there exists a complete  and totally geodesic hypersurface,
orthogonal to $\gamma$  such that  $$u (y)=u( R_\Gamma(y))\quad
\textrm{  for all }\quad  y\in \O.$$  Now since  $u=0$ on  the
boundary $\partial \O$  of $\O$, we deduce from the regularity of $u$  that $\O=R_\Gamma \O$.\\

We are now in position to start the moving "planes" method. From the previous observation,
we  only   need to show that $\O$ is symmetric with respect to a hyperplane orthogonal to the direction $e_1.$ \\

To proceed, we set
\begin{align*}
&D_s:=\{x\in \textrm{H}^{N}, \quad x_1>s\}\\
&\Sigma_s:=D_s\cap \O  \\
&R_s(x):=(2s-x_1, x_2, \dots ,x_N)\quad\textrm{the reflection of  $x$ with respect to}\quad  T_{s}\\
&\Sigma'_s:=R_s \Sigma_s  \\
&d:= \textrm{inf}\{ s\in \R, \quad T_{\mu}\cap
\overline{\O}=\emptyset,\quad  \textrm{for all }\quad \mu >s\}.
\end{align*}
For all  $s\in (0, d)$, we  also consider  the function $$ w_s(x):=
u(R_s(x))-u(x),\quad x\in \Sigma_s.$$ We are going to show  that
there exists a critical position $s_*$ such that
\begin{equation}\label{eq:wvanish}
 w_{s_*}(x)=0\quad \textrm{for all} \quad x\in \Sigma_{s_*}.
\end{equation}

Before we start  proving \eqref{eq:wvanish},  we  make  the
following remark.

\begin{remark}\label{eq: rem2}
Let   $s\in (0, d)$  such that   $\Sigma'_s\subset \O$. We observe
from $(\textrm{H}_1)$  that the function  $u_s: x\mapsto u(R_s(x))$
satisfies the same equation   as $u$  in $\Sigma_s$.  So by
$(\textrm{H}_2)$, the function $ w_s$ satisfies

\begin{equation}\label{inega1}
\mathcal{P}^{-}_{\lambda, \Lambda}(\n_g^2w_s)-k|\n_g w_s|_g-\ell
w_s\leq 0 \quad \textrm{in}\quad  \Sigma_s,
\end{equation}
where $\ell$ is the Lipschitz constant of $f$ on $[0, \max(u)]$.
Furthermore, \eqref{inega1} is invariant under isometry group of the
hyperbolic space.

Let   $x\in \partial \O\cap \partial \Sigma_s$,  we have  $R_s(x)
\in  \O$   since  $\Sigma'_s\subset \O$.  Now recalling that   $u$
is positive on  $\O$  and vanishes on $\partial \O\cap \partial
\Sigma_s$,  it  follows that $ w_s(x)=u(R_s(x))>0.$  Since  $w_s$
vanishes on $T_s$, we have
\begin{align}\label{eq: Linearized}
  \begin{cases}
\mathcal{P}^{-}_{\lambda, \Lambda}(\n_g^2w_s)-k|\n_g w_s|_g-\ell w_s\leq 0 & \quad \textrm{in}\quad  \Sigma_s\vspace{3mm}\\
w_s\geq 0&  \quad\textrm{on }\quad  \partial \Sigma_s.
  \end{cases}
  \end{align}

In addition, combining  Lemma  \ref{Comparapucci}  and  \eqref{eq:
Linearized},  we see that the function  $w_s$ satisfies

\begin{align}\label{eq: Linearizedin euclid}
  \begin{cases}
\mathcal{M}_{\lambda, \Lambda}^{-}(\nabla^2 w_s) -\dfrac{\mu}{x_N}|\nabla w_s|- \dfrac{\ell}{x^2_N} w_s\leq 0 & \quad \textrm{in}\quad  \Sigma_s\vspace{1mm}\\
w_s\geq 0&  \quad\textrm{on }\quad  \partial \Sigma_s.
  \end{cases}
  \end{align}
We also emphasize that the hyperplane $\{x_N=0\}\subset\R^N$ is made
of points at infinity  for the hyperbolic model $(\mathbb{H}^N, g)$,
and since $\O$ is bounded, the coefficients in the first and second
order terms in \eqref{eq: Linearizedin euclid} are bounded.
\end{remark}

We now  begin the proof of  \eqref{eq:wvanish}. \\

We say that the hyperplane $T_s$  has reached the position $s_*< d$
if for all  $\mu \in [s_*, d)$, $w_\mu \geq 0$ in $\Sigma_\mu.$ We
are going to move the hyperplane   $T_s$  varying
continuously  $s$  from the position $d$  to the left. If we prove that the hyperplane  $T_s$   has reached position $s_*$,
then we are done. Since  then we can take  a  hyperplane coming from the other side, that is, starting from $-d$ and moving
to the right. The situation is totally symmetric so the second hyperplane will  also  reach position  $s_*$. This implies $w_{s_*}\geq 0$
 and $w_{s_*} \leq 0$ in $\Sigma_{s_*}$ and hence  $w_{s_*} \equiv0$ in $\Sigma_{s_*}.$\\

We  first show  that the procedure described above can start, that
is, there exists a position $\bar{s}< d$ such that  all  $\mu \in
[\bar{s}, d)$, $w_\mu \geq 0$ in $\Sigma_\mu.$

By \cite[Proposition 2.3]{DaLioandSirakov}, there exists a number
$r>0$ such that for all  $\mu< d$ with $|\Sigma_\mu|<r$, \eqref{eq:
Linearized} satisfies the maximum principle  in $\Sigma_\mu$.

We now fix $\bar{s}< d$ sufficiently close to $d$ such that
$|\Sigma_\mu|<r$ for all  $\mu \in [\bar{s}, d)$. Then  by
\cite[Proposition 2.3]{DaLioandSirakov} and  \eqref{eq: Linearized},
$w_\mu \geq 0$ in $\Sigma_\mu$, for   all  $\mu \in [\bar{s}, d)$.

In addition, we have from  Remark \ref{eq: rem2} that   for all
$s\in (0,d)$,   $w_s>0$ on $\partial \O\cap \partial \Sigma_s$.
Hence by Hopf's lemma  \cite[Proposition 2.6]{DaLioandSirakov}, $w_\mu >0$ in $\Sigma_\mu$, for   all  $\mu \in [\bar{s}, d)$. \\

As a summary, by  reducing continuously  $s$, i.e by moving the
plane  $T_s$  from the position  $d$  to the left, and maintaining
$\Sigma'_s$ inside  $\O$,  the hyperplane $T_s$ will stop at a
critical position $s_*$ defined by
$$s_*:=\inf\{s \leq d, \quad \Sigma'_\mu\subset \O \quad \textrm{and}\quad \langle \nu(z), e_1\rangle> 0\quad
\textrm{for all}\quad \mu >s, \quad z\in T_\mu\cap \partial \O\}$$  where $\nu(z)$ is the outer unit normal at $z\in \partial \O$.
In other terms, the reflection $\Sigma'_s$  of $\Sigma_s$  with
respect to  $T_s$ stays inside $\O$
until one of the following situations occurs:\\

\textbf{ Situation 1}: $\Sigma'_{s_*}$  intersects tangentially  the boundary $\partial \O$  at a point  $P$.\\

\textbf{ Situation 2}:  The hyperplane $T_{s_*}$  becomes orthogonal to   $\partial \O$  at a point  $Q$.\\

We note that $w_s >0$ in $\Sigma_s$ for   all  $s>s_* $, and by
continuity with respect to $s$,  $w_{s_*} \geq 0$  in
$\Sigma_{s_*}$.  Applying the strong maximum principle, it turns out
that

\begin{equation}\label{eqwvanishh}
w_{s_*}\equiv0 \quad \textrm{in} \quad \Sigma_{s_*}\quad
\textrm{or}\quad w_{s_*} >0 \quad \textrm{in} \quad \Sigma_{s_*}.
\end{equation}

We  will  prove that in any of  the   situations above,  $w_{s_*}\equiv0 \quad \textrm{in} \quad \Sigma_{s_*}$.\\

In   situation 1,  we  consider the reflection $P_{s_*}$  of the
point  $P$ with respect to the hyperplane $T_{s_*}$. We have
$w_{s_*}(P_{s_*})= \left.\frac{\partial w_{s_*}}{\partial  \nu}
\right |_{ P_{s_*}}=0. $ Hence by  Hopf's boundary point lemma,
$w_{s_*} \equiv0$ in $\Sigma_{s_*}$.

In the case  of   situation 2, it is no longer  possible  to find a
contradiction  just by analyzing the first and second order
derivatives of $w_{s_*}$  and applying Serrin corner lemma,  valid
for semilinear and quasilinear equations.   Instead, we use  the
fact that  classical solutions to fully nonlinear uniformly elliptic
equations  are $C^{2,\alpha}$ at the boundary and enjoy a Taylor
expansion of order $2+\alpha$, see   \cite[Proposition
2.2]{SilvetreBoyan} and \cite{SilvetreBoyan2}.   The proof of
Theorem \ref{theo2} is achieved applying Proposition \ref{hopfgene}
below,  which provides a  non degeneracy result of order strictly
less than $2+\alpha$ under the assumption  that the ratio
$\Lambda/\lambda$ is close to one.

\section{Proof of Theorem \ref{theo2}}\label{prooftheo2}
The aim of this section is to rule out the corner situation
described above. To proceed, we need  some intermediate results. The
following lemma  is a combination of Proposition  4.2 and Lemma 4.4
in \cite{SilvetreBoyan}.
\begin{lemma}\label{lemma order homoge1}
Let $\sigma\subset \mathbb{S}^{N-1}$  be an open  and  smooth, and
$\mathcal{C}$ be the projected cone
$$\mathcal{C}:=\R^{*}_{+}\cdot\sigma=\{tx: t>0\quad\textrm{and}\quad x\in \sigma\}= \{x\in \R^N\setminus\{0\}:\quad |x|^{-1} x\in \sigma\}.$$
 There exists a number $\beta>0$ and a $\beta$-homogeneous function $\Psi\in C(\overline{\mathcal{C}})$ such that

\begin{equation}\label{eqbeta}
 \mathcal{M}^{-}_{\lambda, \Lambda}(\n^2 \Psi)=0\quad \textrm{and}\quad \Psi>0 \quad \textrm{in}\quad  \mathcal{C},\quad  \Psi=0
 \quad \textrm{on} \quad  \partial \mathcal{C}.\end{equation}
Furthermore, any other solution of this equation is a multiple of
$\Psi$.
\end{lemma}
Note  that the homogeneity constant in Lemma \ref{lemma order
homoge1} is defined (see \cite{SilvetreBoyan}) by
\begin{align}\label{eq:defibeta}
\beta:= \sup \biggl\{ &\widetilde{\beta}>0:  \textrm{there exists a} \quad \widetilde{\beta}  \quad \textrm{homgeneous solution   }
\Phi \in C(\overline{ \mathcal{C}} ) \quad \textrm{of}\nonumber\\
\quad  &\mathcal{M}^{-}_{\lambda, \Lambda}(\n^2 \Phi)\leq0\quad
\textrm{and}\quad \Phi>0 \quad \textrm{in}\quad  \mathcal{C}\biggl
\}.
\end{align}

We also record the following Lemma from   \cite{SilvetreBoyan},
which states that  the homogeneity $\beta   $ of the function $\Psi$
from Lemma \ref{lemma order homoge1} is close to two when
$\mathcal{C}$ is close to the quarter space and the Pucci operator $
\mathcal{M}^{-}_{\lambda, \Lambda}$ is close to the Laplacian.

We consider the quarter space $\Pi$  given  by $$\Pi:=\{x_1>0,
x_N>0\}$$ and set $\pi:=\Pi\cap \mathbb{S}^{N-1}$. Then
$\Pi=\R^{*}_{+}\cdot\pi$ and we have
\begin{lemma}\label{lemma order homoge}
Let $\sigma_n$ be an increasing sequence of smooth  subdomains of
$\pi$ such that $\sigma_n\rightarrow \pi$ as $n\rightarrow\infty$
and $C_n:=\R^{*}_{+}\cdot\sigma_n$. Let   $\Psi_n $  be the
$\beta_n$-homogeneous function given by Lemma \ref{lemma order
homoge1},  applied to the operator  $\mathcal{M}^{-}_{\lambda,
\lambda(1+\frac{1}{n})}$ in  $C_n.$ Then $\beta_n\longrightarrow 2$
as $n \longrightarrow \infty.$
\end{lemma}

For a point $Q\in \R^N$, we define the translated cone of
$\mathcal{C}$ by  $\mathcal{C}(Q):=\mathcal{C}+Q$. The following
result which is fundamental for our contradiction argument is
precisely  \cite[Proposition 4.5]{SilvetreBoyan}, when $Q$ is the
origin in $ \R^N$. See also \cite[Theorem 1.4]{ArmstrongSirakov}.

\begin{proposition}\label{hopfgene}
Let  $\beta$  be  the number in Lemma \ref{lemma order homoge1} for
the cone $\mathcal{C}$. Let $\Sigma \subset \R^{N}$ be a  domain
such that   $Q\in \partial \Sigma$   and   $\Sigma\cap B_{\e_0}(Q)$
is $C^2$ diffeomorphic to $\mathcal{C}(Q)\cap B_{\e_0}(Q)$, for some
$\e_0>0$.  Suppose their exists a nonnegative function  $w \in
C(\overline{\Sigma})$ satisfying
\begin{align}\label{eq:for corner}
 \mathcal{M}^{-}_{\lambda, \Lambda}(\n^2w)-b|\n w|-c w\leq 0& \quad \textrm{in}\quad  \Sigma.
  \end{align}
Then either $w\equiv0$ in $\Sigma$ or  $w>0$ in $\Sigma_0$  and
 \begin{align}\label{eq:hopf at corner}
 \liminf_{ t\searrow0} \frac{w(Q+te)}{t^{\beta}}>0
 \end{align}
 for any direction $e\in \mathbb{S}^{N}$  which enters  $\Sigma.$
\end{proposition}

\proof
A complete proof of Proposition \ref{hopfgene}  is given in \cite[Section 7]{SilvetreBoyan}.
For the reader's  convenience, we  outline the steps needed to recover  \eqref{eq:hopf at corner}.\\

Define the function $\overline{\Psi }$ on $\mathcal{C}(Q)$ by
$\overline{\Psi}(y):=\Psi(y-Q)$, where $\Psi$ is the solution in
Lemma \ref{lemma order homoge1}.  Note that  $\overline{\Psi }$
solves \eqref{eqbeta} in the cone $\mathcal{C}(Q)$ and we have  for
all $t>0$ and all $e\in \mathbb{S}^{N-1}$,

\begin{align}\label{eq:for corner2}
\overline{\Psi}(Q+te)=\Psi(te)=t^{\beta}\Psi(e).
 \end{align}

Next after identifying  $\Sigma$   with  $\mathcal{C}(Q)\cap
B_{\e_0}(Q)$, we define
$$q(r):=\inf_{ \mathcal{C}(Q)\cap (B_{2r}(Q)\setminus B_{r}(Q))}  \frac{w}{\overline{\Psi}},\quad \textrm{for}\quad 0<r<\frac{\e_0}{2}.$$
Following \cite[Section 7]{SilvetreBoyan}, there exists a constant
$C>0$, such that $q(r)\geq C>0$ for  all  small $r$.  This together
with \eqref{eq:for corner2} implies that  for  all small $r$,
$$\inf_{ r\leq t<2r} \frac{w(Q+te)}{t^{\beta}}\geq C>0$$ and \eqref{eq:hopf at corner} follows by letting $r$ goes to zero.
\QED

We are now in position to complete the proof of Theorem \ref{theo2}.

\proof[\textbf{ Proof of Theorem  \ref{theo2} completed}]

We want to show that $$w_{s_*}\equiv0 \quad \textrm{in} \quad
\Sigma_{s_*}.$$ Let us assume in \eqref{eqwvanishh} that   $w_{s_*}
>0 \quad \textrm{in} \quad \Sigma_{s_*}.$ Then from  \eqref{eq:
Linearizedin euclid} and since $\Sigma_{s_*}\subset
(\mathbb{H}^N,g)$,
\begin{align}\label{eq: Linearizedin euclid0}
\mathcal{M}_{\lambda, \Lambda}^{-}(\nabla^2 w_{s_*}) -\mu'|\nabla
w_{s_*} |- \ell' w_s\leq 0 & \quad \textrm{in}\quad  \Sigma_{s_*},
\end{align}
 for some constants $\mu', \ell'>0.$

Because  $\Sigma_{s_*}$ has a right-angle corner at the point $Q$,
then for $\e$ small enough, $\Sigma_{s_*}\cap B_\e(Q)$ is $C^2$
diffeomorphic to a neighborhood of $Q$ in $\Pi$+Q. Hence we can find
a sequence of smooth cones $\mathcal{C}_n$ converging to $\Pi$ from
the inside, and for all $n$, there exists a number  $r_n>0$ such
that
 $\mathcal{C}_n(Q)\cap B_\e(Q)\subset\Sigma_{s_*}\cap B_\e(Q)$ for all $\e\in(0,r_n).$  Applying
 Proposition \ref{hopfgene}  to  \eqref{eq: Linearizedin euclid0} in each $\mathcal{C}_n(Q)\cap B_\e(Q)$, we have
 \begin{equation}\label{eq:woncones}
 w_{s_*}(Q+te)\geq C_nt^{\beta_n}
 \end{equation} for each direction $e$ entering $\mathcal{C}_n(Q)$.
In the other hand  $\n_g^2w_{s_*}(Q)=0=\n_gw_{s_*}(Q)$, see for
instance \cite{Ku-Pra}).  Hence as in  \cite[(3.2)]{SilvetreBoyan},
we apply \cite[Proposition 2.2]{SilvetreBoyan} to \eqref{eq:
Linearizedin euclid0}  and obtain
 \begin{equation}\label{eq:upinequ}
 w_{s_*} (Q+te)\leq C t^{2+\alpha},
 \end{equation}
for every direction $e\in \mathbb{S}^{N}.$ From Lemma \ref{lemma
order homoge},  $\beta_n\rightarrow2$ as $n\rightarrow\infty.$  But
then  \eqref{eq:upinequ} and  \eqref{eq:woncones} are in
contradiction when  $n$ is chosen so large (or $\e$ in Theorem
\ref{theo2} is chosen  so small)  that $\beta_n<2+\alpha.$\QED

\section{Proof of Theorem \ref{theo1}}\label{prooftheo1}
The proof consists  in finding a contradiction in the situation 2
corresponding to the  corner point  in the moving plane method in
Section \ref{Movingplane}.

\proof[\textbf{ Proof of Theorem  \ref{theo1} completed}] Without
loss of general, we assume $s_*=0$.  Let us consider the reflection
$$R(x_1, x_2,\dots, x_N)=(-x_1, x_2,\dots, x_N).$$

For any symmetric matrix   $M=(m_{ij})\in S_N$,  we define
$\overline{M}:=(\e_{ij} m_{ij})$, where $\e_{11}=1,$    $\e_{ij}=1,$
if $i,j\geq2$, $\e_{1j}=-1$ if $j\ne 1$.  We  also put
$\overline{p}:=R(p)$ for any $p\in \R^N.$

Then we  check using \eqref{Hessexpress} that
\begin{equation}\label{eq:forsymme}
\n_g^2u_{s_*}(x)=\overline{\n^2_gu(R(x))}\quad \textrm{and}\quad
\n_g u_{s_*}(x)=\overline{\n_gu(R(x))}  \quad  \textrm{for
every}\quad  x\in \overline{\Sigma_{s_*}}.
\end{equation}
In  particular
\begin{equation}\label{eq:of Dw}
\n_g^2u_{s_*}=\overline{\n^2_gu}\quad \textrm{and}\quad  \n_g
u_{s_*}=\overline{\n_gu}  \quad  \textrm{on }\quad  T_{s_*}\cap
\overline{\Sigma_{s_*}}.
\end{equation}
Since $F$ is continuously differentiable in $M$ and $f$ is locally
Lipschitz, we have that

\begin{equation}\label{eq:of w}
a_{ij}(x)(\n_g^2w)_{ij}+b_{j}(x)(\n_gw)_{j}+c(x)w(x)=0\quad
\textrm{in}\quad \Sigma_{s_*},
\end{equation}
where
\begin{align*}
a_{ij}(x):=&\int_{0}^{1}\frac{\partial F}{\partial m_{ij}} (\n_g^2w^t(x), \n_g w^t(x))\textrm{d}t\\
b_{j}(x):=&\int_{0}^{1}\frac{\partial F}{\partial p_{j}} (\n_g^2w^t(x), \n_g w^t(x))\textrm{d}t\\
c(x):=&\int_{0}^{1}\frac{\partial f}{\partial y} (w^t(x))\textrm{d}t\\
w^t(x):=&tu_{s_*} (x)+(1-t)u(x).
\end{align*}
Now recalling \eqref{eqpuccirelateR}   and \eqref{ChristSymhy},
\eqref{eq:of w} yields
\begin{equation}\label{eq:of w2}
a_{ij}(x)\frac{\partial^2 w}{\partial x_i \partial
x_j}(x)+\overline{b}_{j}(x)\frac{\partial w}{\partial
x_j}(x)+c(x)w(x)=0\quad \textrm{in}\quad \Sigma_{s_*},
\end{equation}
where
\begin{equation}\label{PucciK}
\overline{b}_{j}(x):=
\begin{cases}
2 x_N^{-1} a_{Nj}(x)+ x_N^{2} b_{j}(x)  \quad \textrm{ if} \quad  j=1, \dots, N-1\vspace{3mm}\\
2 x_N^{-1} a_{NN}(x)-\textrm{tr}(A(x))+ x_N^{2} b_{N}(x) \quad
\textrm{ if }  \quad  j=N.
  \end{cases}
\end{equation}
Note  the $a_{ij}$ are continuous since $F$ is $C^1$ in $M$.

Recalling also $(\textrm{H}_1)$, we have
$F(\overline{\n^2_gu(R(x))},
\overline{\n_gu(R(x))})=F(\n^2_gu(R(x)), \n_gu(R(x)))$ and from the
definition of $\overline{M}$,
\begin{equation}\label{eq: aij1}
\frac{\partial F}{\partial m_{1j}}(\overline{\n^2_gu(R(x))},
\overline{\n_gu(R(x))})+\frac{\partial F}{\partial
m_{1j}}(\n^2_gu(R(x)), \n_gu(R(x)))=0, \quad  j>1.
\end{equation}
By change of variable, we have
$$
2a_{ij}(x)=\int_{0}^{1}\frac{\partial F}{\partial m_{ij}}
(\n_g^2w^t(x), \n_g w^t(x))\textrm{d}t+\int_{0}^{1}\frac{\partial
F}{\partial m_{ij}} (\n_g^2\widetilde {w}_t(x), \n_g \tilde
{w}_t(x))\textrm{d}t,
$$
where $\widetilde {w}_t(x):=tu(x)+(1-t)u_{s_*} (x).$ This together
with \eqref{eq:of Dw} and  \eqref{eq: aij1} allow to get
\begin{equation}\label{eq:of a1j}
a_{1j}=0 \quad  \textrm{on }\quad  T_{s_*}\cap
\overline{\Sigma_{s_*}} \quad \textrm{for all }\quad j>1.
\end{equation}
The remaining part of the proof is contained in  \cite[Proof of
Theorem 1.1]{SilvetreBoyan}: using  \eqref{eq:of a1j} and the
continuity of $a_{ij}$, there exists a coordinate system such that
the modified  function, still  denoted by $w_{s_*}$ satisfies
\begin{align}\label{eq:modifiedw}
\mathcal{M}_{\lambda, \Lambda}^{-}(\nabla^2 w_{s_*})
-\tilde{\mu}|\nabla w_{s_*} |- \tilde{\ell} w_s\leq 0 & \quad
\textrm{in}\quad  \widetilde{\Sigma_{s_*}},
\end{align}
with $1-\e<\lambda\leq\Lambda<1+\e$ and $\e>0$. For small  $\e$, the
ration  $\Lambda/\lambda$ is close to one  and we can  repeat the
argument in the proof of  Theorem  \ref{theo2}   to get $w_{s_*}
\equiv 0$ in $\widetilde{\Sigma_{s_*}}$. This complete the proof of
Theorem  \ref{theo1}. \QED

\section{Proof of Theorem \ref{theo1'}}\label{prooftheo1'}
We explain here how to deduce the proof of Theorem  \ref{theo1'}
from the argument  on the hyperbolic space developed in the previous
sections. We  consider the metric $\widetilde{g}_0$ induced on the
unit sphere $\mathbb{S}^N$ from the Euclidean metric   $g_0$ on
$\R^N$. We then transform  problem \eqref{eq:proble1} to an
equivalent problem on $\R^N\setminus\{0\}$ via the stereographic
projection  from  $(\mathbb{S}^N \setminus \mathcal{N},
\widetilde{g}_0)$ to   $\R^N\setminus\{0\}, g)$ for $\mathcal{N}\in
\bar{\O}$, where $g$ is the metric defined by
\begin{equation}\label{eqmetriinflat}
g_{ij}(x)=\frac{4}{(1+|x|^2)^2}\delta_{ij}.
\end{equation}
Proving Theorem \ref{theo1'} then follows  the steps in the proof of Theorem \ref{theo1}
in the previous sections. We also refer the reader to  \cite[Section 3]{Ku-Pra},
where this process is detailed in the particular case  when  the operator $F$ is  the Laplacian of the unit sphere. \\

We emphasize that a key step in the proof of Theorem \ref{theo2}  is
the estimate in Lemma \ref{Comparapucci}. We   also provide the
corresponding estimate  for the metric $g$ in \eqref{eqmetriinflat}.

Using \eqref{ChristSymhy0}, we find
\begin{equation}\label{Chrisy}
\Gamma^{k}_{ij}=-\frac{8}{1+|x|^2}\{x_{i}\delta_{kj}+x_{j}\delta_{ik}-x_{k}\delta_{ij}\},
\end{equation}
and hence
\begin{equation}\label{ChristSymhy}
\nabla^2_g u=\nabla^2 u+\frac{8}{1+|x|^2}x\cdot\n u
\textrm{d}x\otimes \textrm{d}x
-\frac{8}{1+|x|^2}\biggl(\sum^N_{i,j=1}x_{i}\frac{\partial
u}{\partial x_j}\textrm{d}x_i\otimes
\textrm{d}x_j+\sum^N_{i,j=1}x_{j}\frac{\partial u}{\partial
x_i}\textrm{d}x_i\otimes \textrm{d}x_j\biggl).
\end{equation}

We recall that for $a,b\in \R^N$, the symmetric  matrix associated
to $a\otimes b+b\otimes a$  has two nonzero eigenvalues given by
$a\cdot b\pm|a||b|$. Also the Pucci operators  enjoy the properties
$$\mathcal{M}_{\lambda, \Lambda}^{-}(A+B)\geq \mathcal{M}_{\lambda,
\Lambda}^{-}(A)+\mathcal{M}_{\lambda, \Lambda}^{-}(B)\quad
\textrm{and}\quad \mathcal{M}_{\lambda,
\Lambda}^{-}(-A)=\mathcal{M}_{\lambda, \Lambda}^{+}(A)$$ for all
symmetric matrices $A, B$. Using this and $\lambda\leq  \Lambda$, we
find
$$  \mathcal{M}_{\lambda, \Lambda}^{+}\biggl(\sum^N_{i,j=1}x_{i}\frac{\partial u}{\partial x_j}\textrm{d}x_i\otimes \textrm{d}x_j+\sum^N_{i,j=1}x_{j}\frac{\partial u}{\partial x_i}\textrm{d}x_i\otimes \textrm{d}x_j\biggl)\geq 2\lambda x\cdot \n u\geq -\lambda |x||\n u|$$ which implies
$$\mathcal{M}_{\lambda, \Lambda}^{-}(\nabla^2_g u)\geq \mathcal{M}_{\lambda, \Lambda}^{-}(\nabla^2 u)-8N\lambda\frac{|x|}{1+|x|^2}|\n u|-8\lambda\frac{|x|}{1+|x|^2}|\n u|. $$

Finally, it follows  from \eqref{relagen} and  \eqref{eqmetriinflat}
that
\begin{equation}\label{relapuccii}
 \mathcal{P}_{\lambda, \Lambda}^{-}(\nabla^2_g u)=\frac{(1+|x|^2)^2}{4} \mathcal{M}_{\lambda, \Lambda}^{-}(\nabla^2_g u)\geq  \frac{(1+|x|^2)^2}{4} \mathcal{M}_{\lambda, \Lambda}^{-}(\nabla^2 u)-2\lambda (N+1)|x|(1+|x|^2)|\n u|.
\end{equation}
The coefficients in  \eqref{relapuccii}  are   smooth  in
$\R^N\setminus\{0\}$ and  $\nabla^2_g u$  satisfies
\eqref{eq:forsymme}. Therefore the  procedure in Sections
\ref{Movingplane} and \ref{prooftheo1} remains valid in
$(\R^N\setminus\{0\}, g)$.


\begin{thebibliography}
\footnotesize

\bibitem{Alexandrov} A. D. Alexandrov,  {Uniqueness Theorem
for surfaces in large I}, Vestnik Leningrad Univ. Math. 11 (1956),
5-17.

\bibitem{ArmstrongSirakov}   S.N. Armstrong,  B. Sirakov and  C.K. Smart,  Singular solutions of fully nonlinear elliptic equations and
applications. Arch. Ration. Mech. Anal. 205(2); 345-394 (2012).


\bibitem{BCN} H. Berestycki, L. A. Caffarelli and L. Nirenberg, Monotonicity for elliptic equations in unbounded
Lipschitz domains. Comm. Pure Appl. Math. 50 (1997) 1089-1111.

\bibitem{Berestycki and L} H. Berestycki and L. Nirenberg, Monotonicity, symmetry and antisymmetry of solutions of semilinear
elliptic equations. J. Geom. Phys., 5(2):237-275, 1988.

\bibitem{Berestycki}  H. Berestycki and L. Nirenberg. Some qualitative properties of solutions of semilinear elliptic equations in cylindrical domains. In Analysis, et cetera: Research papers published in honor of Jrgen
Moser's 60th birthday, pages 255-273. Academic Press, 1990.

\bibitem{BerestC} H. Berestycki and L. Nirenberg. On the method of moving planes and the sliding method. Bol. Soc.
Brazil (N.S.),  22(1):1-37, 1991.

\bibitem{DemengelBir} I. Birindelli and F. Demengel. Overdetermined problems for some fully non linear operators. Comm. Part. Diff. Eq., 38(4):608-628, 2013.

\bibitem{Brandolini} B. Brandolini, C. Nitsch, P. Salani, C. Trombetti, Serrin type overdetermined problems: an alternative proof. Arch. Rational Mech. Anal. 190 (2008), 267-280.

\bibitem{Birindeli-Demengel} Birindelli, I., and F. Demengel. Overdetermined problems for some fully non linear operators. Communications in Partial Differential Equations 38.4 (2013): 608-628.

\bibitem{CiraoloLuigi} G. Ciraolo, L. Vezzoni, A rigidity problem on the round sphere, Commun. Contemp.
Math. 19, no. 05 (2017): 1750001.

\bibitem{DamascelliA} L. Damascelli and F. Pacella. Monotonicity and symmetry of solutions of p-Laplace equations,
$1 < p < 2$, via the moving plane method. Ann. Scuola Norm. Sup.
Pisa (4), 26(2):689-707, 1998.

\bibitem{DamascelliB} L. Damascelli and B. Sciunzi. regularity, monotonicity and symmetry of positive solutions of m-Laplace equations. J. Differential Equations, 206(2):483-515, 2004.

\bibitem{DancerS} E. N. Dancer. Some notes on the method of moving planes. Bull. Austral. Math. Soc., 46(3):425-434,
1992.

\bibitem{DS} E. Delay and P. Sicbaldi, Extremal domains for the first
eigenvalue of the Laplace Beltrami operator in a general compact
riemannian manifold,  In Annales de l'institut Fourier, vol. 59,
no.2, pp. 515-542. 2009.



\bibitem{DaLioandSirakov}  F. Da Lio and B. Sirakov, Symmetry results for viscosity solutions of fully nonlinear uniformly elliptic equation.
 Journal of the European Mathematical Society 9, no. 2 (2007): 317-330.


\bibitem{VazquezEncicoSalas} M. Dominguez-Vazquez, A. Enciso, and  D. Peralta-Salas, Solutions to the overdetermined boundary problem for semilinear equations with position-dependent nonlinearities. Advances in Mathematics 351 (2019): 718-760.


\bibitem{FallIgnace} M. M. Fall and I. A. Minlend, \textsl{Serrin's
over-determined problem in Riemannian manifolds.\/} Adv. Calc. Var.
{\bf 8} (2015), 371-400.


\bibitem{Fall-MinlendI-Weth}  M.M. Fall, I. A. Minlend, T. Weth, Unbounded periodic solutions to Serrin's overdetermined boundary value problem. Arch. Ration. Mech. Anal. 233(2017), no. 2, 737-759.

\bibitem{FallMinlendIWeth2}  M.M. Fall, I.A. Minlend, T. Weth, Serrin overdetermined problem on the sphere. Cal. Var. Partial Differential Equations 57, no. 1 (2018): 3.

\bibitem{farina-valdinoci:2013-2} A. Farina, L. Mari,  E. Valdinoci, Splitting theorems, symmetry
results and overdetermined problems for Riemannian manifolds. Comm.
Partial Differential Equations 38 (2013), no. 10, 1818-1862.

\bibitem{Garofalo} N. Garofalo, J.L. Lewis, A symmetry result related to some overdetermined boundary
value problems, Amer. J. Math. 111 (1989), 9-33.

\bibitem{BGidasW}
B. Gidas, W.M. Ni, and L. Nirenberg. Symmetry and related properties
via the maximum principle. Comm. Math. Phys., 68(3):209-243, 1979.

\bibitem{Gidas}  B. Gidas, W.M. Ni, and L. Nirenberg. Symmetry of positive solutions of nonlinear elliptic equations in $\R^n$. In Mathematical analysis and applications, Part A, volume 7 of Adv. in Math. Suppl. Stud., pages 369-402. Academic Press, 1981.


\bibitem{Ku-Pra}  S. Kumaresan, J. Prajapat, Serrin's result for hyperbolic space and sphere. Duke Math. J. 91 (1998), 17-28.

 \bibitem{Molzon} R. Molzon, Symmetry and overdetermined boundary value problems. Forum Math. 3 (1991), 143-156.

\bibitem{FPP} F. Pacard and P. Sicbaldi,  Extremal domains for the first eigenvalue of the Laplace-Beltrami operator. Ann. Inst. Fourier (Grenoble)
59 (2009), no. 2, 515-542.


\bibitem{JPau}  J.  Paupert, Introduction to Hyperbolic Geometry.  Arizona State University Lecture Notes (2016).


\bibitem{Poho} S. I. Pohozaev, On the eigenfunctions of the equation $\Delta u+f(u) = 0,$ Dokl. Akad. NaukSSSR 165 (1965), 36-39.


\bibitem{PucciSerrin} P. Pucci and J. Serrin, The maximum principle. Progress in Nonlinear Differential Equations and Their Applications, Birkhauser, Basel, 2007.

\bibitem{iuXia} G. Qiu and C. Xia. (2017), Overdetermined boundary value problems in Sn. J. Math. Study, 50(2), 165-173.


\bibitem{Ros-Ruiz-Sicbaldi-2015} A. Ros, D. Ruiz, P. Sicbaldi,
A rigidity result for overdetermined elliptic problems in the plane.
Communications on Pure and Applied Mathematics 70.7 (2017):

\bibitem{SilvetreBoyan}  L. Silvestre  and B. Sirakov,  Overdetermined problems for fully nonlinear elliptic equations.  Calculus of Variations and Partial Differential Equations 54.1 (2015): 989-1007.

\bibitem{SilvetreBoyan2}  L. Silvestre  and B. Sirakov,  Boundary regularity for viscosity solutions of fully nonlinear elliptic equations.  Communications in Partial Differential Equations 39.9 (2014): 1694-1717



1223-1252.

\bibitem{Serrin} J. Serrin, A Symmetry Theorem in Potential Theory. Arch. Rational Mech. Anal. 43 (1971), 304-318.

\bibitem{Sic} P. Sicbaldi,  New extremal domains for the first eigenvalue
of the Laplacian in flat tori. Calc. Var. (2010) 37:329-344.

\bibitem{Wein} H. F. Weinberger, Remark on the preceding paper of Serrin. Arch. Rational Mech. Anal. 43 (1971), 319-320.
\end{thebibliography}
\end{document}